%Simbolos matematicos.

%Macros.

\def\ouv#1{\smash{\mathop{#1}\limits^{\lower 1pt\hbox
{$\scriptscriptstyle\circ$}}}}

\def\hfl#1#2{\smash{\mathop{\hbox to 12mm{\rightarrowfill}}
\limits^{\scriptstyle#1}_{\scriptstyle#2}}}

%Titulos, enunciados.

\def\stit#1{\vskip 3mm plus 1mm minus 2mm {\bf{#1}}
		\smallskip}

\font\tir=cmbx10 at 12pt

\def\ref#1#2#3#4{{\bf #1}{\ #2}{\it ,\ #3}{,\ #4}\medskip}

%Dibujos

\def \picture #1 by #2 (#3){\midinsert \centerline 
{\vbox to #2{\hrule width #1 heigth 0pt 
depth 0pt \null \vfill \special {picture #3}}}\endinsert }

\def\scaledpicture #1 by #2 (#3 scaled #4) {{
\dimen0 =#1 \dimen1 =$2
\divide \dimen0 by 1000 \multiply \dimen0 by #4
\divide \dimen1 by 1000 \multiply \dimen1 by #4
\picture \dimen0 by \dimen1 (#3 scaled $4)}}

\def\figure #1 #2 #3 {\midinsert \vglue 3mm 
{\vbox to #3 {\hrule width 6cm height 0cm depth 0cm \vfill
{\special {picture #1 scaled #2}}}}\vglue 2mm \endinsert}

\magnification=1200

\input psfig.sty

\overfullrule=0pt

{\centerline {\tir {Maximal Abelian Torsion Subgroups of Diff({\bf
C},0) }}}

\bigskip
\bigskip

{\centerline {Kingshook Biswas \footnote{$^1$} {Ramakrishna Mission Vivekananda University,
Belur Math, WB-711202, India.
E-mail: kingshook@rkmvu.ac.in} }}

\bigskip
\bigskip

{\bf Abstract.} {\it In the study of the local dynamics of a germ
of diffeomorphism fixing the origin in {\bf C}, an important
problem is to determine the centralizer of the germ in the group
Diff({\bf C},0) of germs of diffeomorphisms fixing the origin.
When the germ is not of finite order, then the centralizer is
abelian, and hence a maximal abelian subgroup of Diff({\bf C},0).
Conversely any maximal abelian subgroup which contains an element
of infinite order is equal to the centralizer of that element. A
natural question is whether every maximal abelian subgroup
contains an element of infinite order, or whether there exist
maximal abelian torsion subgroups; we show that such subgroups do
indeed exist, and moreover that any infinite subgroup of the
rationals modulo the integers {\bf Q/Z} can be embedded into
Diff({\bf C},0) as such a subgroup.}

\medskip

{\centerline {\it AMS Subject Classification: 37F99}}

\bigskip

\bigskip

\stit {1. Introduction.}

\bigskip

We consider the group of germs of diffeomorphisms in {\bf C}
fixing the origin, Diff({\bf C},$0) = \{ f : f(z) = \lambda z +
O(z^2), \lambda \neq 0\}$. The local dynamics of such germs $f$
near the fixed point $0$ has been intensively studied, in particular
the question, when is $f$ {\it linearizable}, i.e. conjugate to its linear
part $L_{\lambda}(z) = \lambda z$. If the fixed point is {\it attracting}
($|\lambda| < 1$) or {\it repelling} ($|\lambda| > 1$) then a classical theorem
of Koenigs asserts that $f$ is linearizable. When the fixed point is {\it indifferent}
($|\lambda| = 1, \lambda = e^{2\pi i \alpha},
\alpha \in ${\bf (R/Z}), the linearizability of $f$ depends very sensitively on
the arithmetic of the {\it rotation number} $\alpha$. Any {\it nondegenerate parabolic}
germ ($\alpha = p/q \in ${\bf Q/Z}, $f^q \neq id$) is nonlinearizable whereas a
{\it degenerate parabolic} germ ($\alpha = p/q \in ${\bf Q/Z}, $f^q = id$) is always linearizable.
When $\alpha$ is irrational, for $\alpha$ poorly approximable by
rationals all germs $f$ are linearizable whereas for $\alpha$
very well approximable by rationals there exist nonlinearizable
germs. The sharp arithmetic condition is called the Brjuno condition
([Br]); its optimality was shown by Yoccoz ([Yo]) (see for example [PM1] for a survey of the linearization
problem).

\medskip

The centralizer Cent($f) = \{ g : g
\circ f = f \circ g \}$ of a germ $f$ in Diff({\bf C},0) can be thought
of as the group of symmetries of the dynamics (it's elements
conjugate the dynamics to itself). The centralizer clearly contains
any abelian subgroup containing $f$; when $f$ is not of finite
order, it is well known that Cent($f$) is abelian (we recall the
description of such centralizers in the following section), and is
hence a maximal abelian subgroup of Diff({\bf C},0). Moreover when $f'(0) = \lambda$ is not a
root of unity then the group homomorphism given by the ''rotation number map''
$$\eqalign{
\rho :  \hbox{Diff}({\bf C},0) & \to {\bf C/Z} \cr
                g              &\mapsto {1 \over 2\pi i} \log g'(0) \cr
}$$
is injective restricted to Cent($f$), identifying the centralizer
with a subgroup of {\bf C/Z}. The restriction to Cent($f$) is
surjective if and only if $f$ is linearizable, whereas if $f$ is
nonlinearizable, then  $\rho($Cent($f)) \subset ${\bf R/Z}. Understanding
the arithmetic of subgroups of {\bf R/Z} which occur as groups of rotation numbers
$\rho($Cent($f))$ in the nonlinearizable case can thus be seen as a
generalization of the linearization problem. This seems to be a
very difficult problem for which few results are known: Moser [Mo]
has shown that the irrationals occurring in such subgroups must
admit good simultaneous rational approximations, while Perez-Marco
[PM2] has constructed examples where the subgroups are
uncountable, containing Cantor sets.

\medskip

Any maximal abelian subgroup of Diff({\bf C},0)
containing an infinite order element $f$ is equal to Cent($f$). In
order to classify all maximal abelian subgroups of Diff({\bf
C},0), it is natural to ask therefore whether there exist maximal
abelian torsion subgroups (otherwise only centralizers can occur).
We show that this is indeed the case, and moreover any infinite
subgroup of {\bf Q/Z} can occur as the corresponding group of rotation numbers:

\medskip

\noindent {\bf Theorem 1.}{\it For any infinite subgroup $H$ of
{\bf Q/Z}, there exists a maximal abelian torsion subgroup
$\hat{H}$ of Diff({\bf C},0) such that $\rho$ restricted to
$\hat{H}$ is injective, and $\rho(\hat{H}) = H$.}

\medskip

Thus the rotation numbers of maximal abelian torsion subgroups
can be arbitrary and are not subject to any arithmetic condition.

\bigskip

\bigskip

\stit {2. Centralizers in Diff({\bf C},0).}

\bigskip

We recall some well known facts about centralizers of elements of
Diff({\bf C},$0$). We denote by {\bf C}[[$z$]] the ring of formal
power series in $z$ and by Diff$_{For}$({\bf C},0) $:= \{ \ f \in \hbox{\bf C}[[z]] :
f(0) = 0, f'(0) \neq 0 \ \}$ the group of formal diffeomorphisms of germs fixing
$0$. We identify Diff({\bf C},0) is with the subgroup of elements
of Diff$_{For}$({\bf C},0) whose series converge. For $f \in$
Diff({\bf C},0) we denote its centralizer in Diff$_{For}$({\bf C},0)
by Cent$_{For}(f)$; the analytic centralizer Cent$(f)$ is then
identified with the subgroup of elements of the formal centralizer
Cent$_{For}(f)$ whose series converge.

\medskip

\noindent {\bf Proposition 2.}{\it For any $\lambda \in ${\bf C} $^*$
 which is not a root of unity, Cent$(L_{\lambda}) =
$Cent$_{For}(L_{\lambda}) = \{ \ L_{\mu} \}_{\mu \in {\bf C}^*}$.}

\medskip

\noindent {\bf Proof:} Comparing coefficients of both sides of $g
\circ L_{\lambda} = L_{\lambda} \circ g$ gives $\lambda^n g_n =
\lambda g_n, n \geq 1$ (where $g(z) = \sum_n g_n z^n$),
so $g_n = 0$ for $n \geq 2$. $\diamond$

\medskip

We recall that when $\lambda$ is not a root of unity, any $f(z) =
\lambda z + f_2 z^2 + \dots \in$Diff$_{For}$({\bf C},0) is
formally linearizable; there exists a unique formal germ $h_f(z) =
z + h_2 z^2 + \dots \in$Diff$_{For}$({\bf C},0) conjugating $f$ to
$L_{\lambda}$, and $f$ is linearizable if and only if the formal series $h_f$
is convergent. Indeed comparing coefficients of both sides of the
conjugacy equation $h_f \circ f = L_{\lambda} \circ h_f$
determines a recursive solution of the form $h_n = {P_n(\lambda,
f_2, \dots, f_n, h_1, \dots h_{n-1}) \over (\lambda^n - \lambda)}, n \geq
2$, where the $P_n$'s are universal polynomials.

\medskip

Observing that any germ (formal or analytic) conjugating two germs
also conjugates their centralizers, it follows from the previous
proposition that

\medskip

\noindent{\bf Proposition 3.}{\it \ For any $f(z) = \lambda z + O(z^2)
\in$ Diff({\bf C},0) with $\lambda$ not a root of unity, Cent$_{For}(f)
= \{ h_f \circ L_{\mu} \circ h^{-1}_f : \mu \in
{\bf C}^* \}$. If

\noindent (i) $f$ is linearizable then Cent($f) =$ Cent$_{For}(f) = \{ h_f \circ L_{\mu} \circ h^{-1}_f : \mu \in
{\bf C}^* \}$ and $\rho : $Cent$(f) \to {\bf C/Z}$ is an isomorphism.

\noindent (ii) $f$ is nonlinearizable then Cent($f) = \{ h_f \circ L_{\mu} \circ h^{-1}_f
: \mu \in {\bf C}^* \ \hbox{such that} \ h_f \circ L_{\mu} \circ h^{-1}_f \ \hbox{converges}
\}$, all germs in Cent$(f)$ are indifferent and $\rho : $Cent$(f) \to {\bf R/Z}$ is an
injective homomorphism.}

\medskip

In either case above we see that the centralizer is abelian (and
hence a maximal abelian subgroup of Diff({\bf C},0)). In the case
of a nondegenerate parabolic germ $f$, while $f$ is not (even
formally) linearizable, it can always be formally
embedded into the flow of a holomorphic vector field with
a zero at the origin. Indeed for $n \geq 1, \tau \in ${\bf C},
the vector field $X_{n,\tau} = {z^{n+1} \over 1 + \tau z^n} {\partial \over
\partial z}$ is a normal form for any holomorphic vector field with a zero of order
(n+1) at the origin, and if $\rho(f) = p/q, f^q(z) = z + az^{n+1} + \dots$ with
$a \neq 0$, then for some $\tau = \tau(f)$ there exists a formal germ $\phi$ such that
$f = \phi \circ e^{2\pi i p/q} \hbox{exp}(X_{n,\tau}) \circ
\phi^{-1}$. The formal and analytic centralizers of exp$(X_{n,\tau})$
are equal, both given by the abelian group of germs $\{ \ e^{2\pi i
k/n} \hbox{exp}(tX_{n,\tau}) : k \in {\bf Z}/n{\bf Z} , t \in {\bf
C} \ \} \cong {\bf Z}/n{\bf Z} \times {\bf C}$ (note the rotations
around the origin of finite order $n$ commute with the flow of $X_{n,\tau}$).
For these classical results we refer the reader to the articles of Baker ([Ba]),
Ecalle ([Ec]) and Voronin ([Vo]). They allow us to describe the
centralizer of nondegenerate parabolic germs as follows:

\medskip

\noindent {\bf Proposition 4.} {\it \ For any nondegenerate
parabolic germ $f$ such that $\rho(f) = p/q, f^q(z) = z + az^{n+1} + \dots, a \neq
0$, for some $\tau \in${\bf C}, Cent$(f)$ is given by the subgroup of elements of  $\{ \ \phi \circ e^{2\pi i
k/n} \hbox{exp}(tX_{n,\tau}) \circ \phi^{-1} : k \in {\bf Z}/n{\bf Z} , t \in {\bf
C} \ \}$ which converge (where $\phi \in $Diff$_{For}$({\bf C},0)
formally conjugates $f$ to $e^{2\pi i p/q}
\hbox{exp}(X_{n,\tau})$).}

\medskip

Thus the centralizer is again abelian (isomorphic to a subgroup of
${\bf Z}/n{\bf Z} \times {\bf C}$), hence maximal abelian, all
elements of the centralizer are parabolic, and $\rho(Cent(f))
\subseteq \left({1 \over n}{\bf Z} \right)/ {\bf Z}$ is finite.

\bigskip

\bigskip

\stit {3. Maximal Abelian Torsion Subgroups of Diff({\bf C},0).}

\bigskip

For $r>0$ we denote by {\bf D}$_r$ the disc of radius $r$ centered
around the origin. For any simply connected domain $D$ containing
$0$ and $\alpha \in ${\bf R/Z}, we denote by $R_{D,\alpha}$ the
intrinsic rotation of $D$ around $0$ by angle $2\pi \alpha$, i.e.
the unique automorphism $R:D \to D$ such that $R(0) = 0, R'(0)=
e^{2\pi i \alpha}$, which can be described as the conjugate of the rigid rotation
$R_{\alpha}(z) = e^{2\pi i \alpha}z$ by any Riemann mapping $h : {\bf D} \to D$
such that $h(0) = 0$. The intrinsic rotations of a given domain
form a group isomorphic to {\bf R/Z}, and any conformal mapping between
two simply connected domains which fixes the origin conjugates
their intrinsic rotations. The construction of maximal abelian
torsion subgroups rests on the following

\medskip

\noindent{\bf Proposition 5.} {\it Let $r, M, \delta >0$ be real and $q, a\geq
2$ integers. Given $f \in $Diff({\bf C},0) such that $\rho(f) = 1/q, f^q = id$,
there exists $\phi \in $Diff({\bf C},0), $\phi = \phi(r, M, \delta, a, f)$
such that $\rho(\phi) = 1/(aq), f = \phi^a$, and for any $g \in $Diff({\bf C},0),
if $g,g^{-1}$ are univalent on ${\bf D}_r$, $|g'(0)|, |(g^{-1})'(0)| \leq M$, and
$g$ commutes with $\phi$, then $\rho(g) \in \left({1 \over q}{\bf Z}+{\bf D}_{\delta}\right)/{\bf
Z}$ (where ${1 \over q}${\bf Z}$+{\bf D}_{\delta}$ denotes numbers of the form
$k/q + \tau, k \in ${\bf Z}, $|\tau| < \delta$).}

\medskip

\noindent{\bf Proof:} The idea is roughly, given $f$, to find an $\phi$ satisfying
$f = \phi^a$ and having singularities very close to the
origin which are almost symmetrically distributed with respect to
the rotation $R_{1 \over q}$, so that any germ commuting with
$\phi$ must preserve the singularities and hence have
rotation number close to a multiple of $1/q$. We achieve this as
follows:

\medskip

Fix a linearization $h(z) = z + O(z^2)$ of $f$, so $f = h^{-1} \circ R_{1/q} \circ h$.
Let $\epsilon_0 , \epsilon_1 > 0$ be such that $h^{-1}$ is univalent
on ${\bf D}_{\epsilon_0}$ and ${\bf D}_{\epsilon_1} \subset U = h^{-1}({\bf D}_{\epsilon_0})
\subset {\bf D}_r$. For $0 < \epsilon < \epsilon_0$
let $D(\epsilon, q)$ be the slit domain $D(\epsilon,
q) := {\bf D}_{\epsilon_0} - \cup_{j \in {\bf Z}/q{\bf Z}} \{ t \, e^{2\pi i
j/q} : \epsilon \leq t < \epsilon_0 \}$. Note $R_{1/q}(D(\epsilon,
q)) = D(\epsilon,q)$, so $R_{D(\epsilon,q), 1/q} = R_{1/q} = h \circ f \circ h^{-1}$; thus
if we set $\phi := h^{-1} \circ R_{D(\epsilon,q), 1/(aq)} \circ h$
then $f = \phi^a$. We check that for $\epsilon$ small enough,
depending only on $r,\delta,a$ and $h$, the germ $\phi$ satisfies
the conclusion asserted by the Proposition:

\medskip

The domain $V_{\epsilon} := h^{-1}(D(\epsilon, q))$
is a slit domain equal to $U$ minus a finite union of slits
$(\gamma_j)_{j \in {\bf Z}/q{\bf Z}}$, which is invariant under
$\phi$. Let $z_j = h^{-1}(\epsilon \, e^{2\pi i j/q})$ be the endpoints
of the slits $\gamma_j$. Then any $z^* \in \gamma_j$ distinct from
$z_j$ is biaccessible from $V_{\epsilon}$, and there are two paths
$\beta, \beta'$ in $V_{\epsilon}$ landing at $z^*$ such that any
conformal representation from $V_{\epsilon}$ to {\bf D} tends to distinct
points of $\partial {\bf D}$ as $z$ tends to $z^*$ along $\beta,
\beta'$. In the plane of $w = h(z)$, the corresponding point
$w^* = h(z^*)$ is biaccessible from $D(\epsilon, q)$, and any intrinsic rotation
$R_{D(\epsilon,q), \alpha}$ with $\alpha \notin {1 \over q}{\bf
Z}/{\bf Z}$ tends to distinct limits in $\partial D(\epsilon,q)$
as $w=h(z)$ tends to $w^* = h(z^*)$
along $h(\beta), h(\beta')$. So the conjugate germ $\phi$ will
tend to distinct limits contained in $\partial V_{\epsilon} \subset {\bf D}_r$ as $z$ tends
to $z^*$ along $\beta, \beta'$.

\medskip

Now let $g(z) = \mu z + O(z^2)$ be a germ such that $g, g^{-1}$
are univalent in ${\bf D}_r$ and $|g'(0)|, |(g^{-1})'(0)| \leq M$. By
classical results on univalent functions, the family ${\cal
F}_{M,r}$ of such functions is a normal family, so for $\epsilon$
small enough depending only on $r,M$, $g({\bf D}_{2\epsilon})
\subset {\bf D}_{\epsilon_1} \subset U$. Suppose $g$ commutes with
$\phi$. Then taking $z^* \in \gamma_j \cap {\bf D}_{2\epsilon}, z^* \neq
z_j$, and letting $z$ tend to $z^*$ along
$\beta,\beta'$ in the equation $\phi(g(z)) = g(\phi(z))$, since
the RHS tends to two distinct limits, we must have $g(z^*) \in \gamma_{j'}$
for some $j'$ (otherwise the LHS would tend to a unique limit). It
follows that $g({\bf D}_{2\epsilon} \cap (\cup_j \gamma_j))
\subseteq g({\bf D}_{2\epsilon}) \cap (\cup_j \gamma_j)$,
applying the same argument to $g^{-1}$ gives the
reverse inclusion, so $g({\bf D}_{2\epsilon} \cap (\cup_j \gamma_j)) =
g({\bf D}_{2\epsilon}) \cap (\cup_j \gamma_j)$. In particular, for any
$j$, $g(z_j) = z_{j'}$ for some $j'$. So
$$
{g(z_j) \over z_j} = {z_{j'} \over z_j} = {h(\epsilon \, e^{2\pi i j'/q}) \over h(\epsilon \, e^{2\pi i j/q})}
$$
Since $g$ belongs to the normal family ${\cal
F}_{M,r}$, we have a uniform estimate  $|g(z_j)/z_j - \mu| \leq C|z_j|$
where the constant $C$ only depends on $r,M$ and not on $g$.
It is clear that $|z_j| = O(\epsilon)$ and
that the RHS is equal to $e^{2\pi i (j' - j)/q} + O(\epsilon)$
with the constants in the $O(\epsilon)$ terms only depending on
$h$; it follows that
$$
\mu = e^{2\pi i (j' - j)/q} + O(\epsilon)
$$
with the constant in the error term $O(\epsilon)$ only depending
on $r,M$ and $h$, so for $\epsilon$ small enough, depending only
on these parameters and $\delta$ but not on $g$, we will have
$$
\rho(g) = {1 \over 2\pi i} \log \mu \in \left({(j'-j) \over q}+{\bf D}_{\delta}\right)/{\bf
Z} \qquad \qquad \diamond.
$$

\medskip

We need the following simple description of subgroups of {\bf
Q/Z}:

\medskip

\noindent{\bf Proposition 6.}{\it (i) Any finite nonzero subgroup $H$ of {\bf
Q/Z} is cyclic, $H = <1/q>$ for some $q \geq 2$.

\noindent (ii) Any infinite subgroup $H$ of {\bf
Q/Z} is an increasing union of cyclic subgroups, $H = \cup_n
<1/q_n>$ for some sequence $(q_n)$ such that $q_n | q_{n+1}$ for
all $n$.}

\medskip

\noindent{\bf Proof:} (i) Let $x_0 = p/q \in H, (p,q$ coprime, $q\geq
2$) be such that $d(x_0, {\bf Z}) = $Min$_{x \in H, x \neq 0} d(x, {\bf
Z})$. Then since $<p/q> = <1/q>$ in {\bf Q/Z}, we must have $x_0 =
\pm 1/q$. For any $x = p'/q' \in H (p',q'$ coprime, $q'\geq
2$), $<1/q'> = <p'/q'> \subseteq H$ and we can write $1/q' = a/q +
r$ for some integer $a$ and some $0 \leq r < 1/q$. Then $r = 1/q'
- a/q \in H$ and $d(r,{\bf Z}) = r < 1/q = d(x_0,{\bf Z})$ so $r =
0$, hence $1/q' = a/q$ and $x \in <x_0>$. Thus $H = <1/q>$.

\medskip

(ii) Let $(x_n)_{n\geq 1}$ be an enumeration of $H$. Then $H$ is
the increasing union of the finite subgroups $H_n = <x_1, \dots,
x_n>$, each of which is of the form $H_n = <1/q_n>$ by (i), and
$q_n = O(H_n)$ divides $q_{n+1} = O(H_{n+1})$ since $H_n$ is a
subgroup of $H_{n+1}$. $\diamond$

\medskip

We can now construct maximal abelian torsion subgroups as follows:

\medskip

\noindent{\bf Proof of Theorem 1:} Given an infinite subgroup $H$
of {\bf Q/Z}, write it in the form $H = \cup_{n \geq 1} <1/q_n>$ where
$q_{n+1} = a_{n+1} q_n$ for some integers $a_n \geq 2$ (we may
assume the cyclic groups are strictly increasing).

\medskip

Fix a monotone decreasing sequence $(r_n)$ converging to $0$,
and an increasing sequence $M_n$ converging to $+\infty$. Let
$\delta_n = {1 \over 3}{1 \over q_{n+1}}$. It is easy to check
that for this choice of $\delta_n$,
$$
\bigcap_{n \geq N} \left( {1 \over q_n}{\bf Z} + {\bf D}_{\delta_n}
\right)/{\bf Z} = \left( {1 \over q_N}{\bf Z} \right)/{\bf Z}
$$
for any $N \geq 1$.

\medskip

Let $f_1 = R_{1/q_1}$ and for $n \geq 1$ define $f_{n+1}$ inductively
by $f_{n+1} = \phi(r_n, M_n, \delta_n, a_{n+1}, f_n)$ (where $\phi$ is
the germ given by Proposition 5). Then $\rho(f_n) = 1/q_n, f_n = f^{a_{n+1}}_{n+1}$,
so the increasing union of cyclic subgroups  $\hat{H} := \cup_{n \geq 1}
<f_n>$ is an abelian torsion subgroup of Diff({\bf C},0) such that
$\rho(\hat{H}) = H$. We check that $\hat{H}$ is maximal abelian:

\medskip

Let $g$ be a germ commuting with $f_n$ for all $n$. Let $N = N(g) \geq 1$
be such that $g$ and $g^{-1}$ are univalent
on {\bf D}$_{r_n}$ and $|g'(0)|, |(g^{-1})'(0)| \leq M_n$ for all $n \geq N$. By the choice of the
$f_n$'s and Proposition 5, it follows that
$$
\rho(g) \in \bigcap_{n \geq N} \left( {1 \over q_n}{\bf Z} + {\bf D}_{\delta_n}
\right)/{\bf Z} = \left( {1 \over q_N}{\bf Z} \right)/{\bf Z}
$$
Thus $g$ is parabolic. Since $\rho($Cent$(g)) \supseteq
\rho(\hat{H}) = H$ is infinite, $g$ must be a degenerate
parabolic, $g^{q_N} = id$. If $\rho(g) = k/q_N$, then the germ $g
\circ f^{-k}_{q_N}$ is tangent to the identity, of finite order
(since $g, f_{q_N}$ commute and are of finite order), and hence
must be the identity. So $g = f^k_{q_N} \in \hat{H}$. $\diamond$

\medskip

To summarize, we can classify maximal abelian subgroups of
Diff({\bf C},0) according to their groups of rotation numbers as
follows:

\medskip

\noindent {\bf Theorem 7.}{\it Any maximal abelian subgroup $\hat{H}$ of
Diff({\bf C},0) must be of one of the following kinds (and all
kinds occur): either

\noindent (i) $\rho(\hat{H})$ is a finite subgroup of {\bf Q/Z}
in which case $\hat{H}$ contains a nondegenerate parabolic germ $f$, is equal
to Cent$(f)$, and is isomorphic to a subgroup of ${\bf Z}/n{\bf Z} \times {\bf
C}$,

or

\noindent (ii) $\rho(\hat{H})$ is an infinite subgroup of {\bf
Q/Z} in which case all elements of $\hat{H}$ are degenerate
parabolic germs, and $\hat{H}$ is isomorphic to the subgroup
$\rho(\hat{H})$ of {\bf Q/Z},

or

\noindent (iii) $\rho(\hat{H}) \cap {\bf (R - Q)/Z} \neq
\emptyset$, in which case $\hat{H}$ contains an irrationally indifferent
 germ $f$, is equal to Cent$(f)$ and isomorphic to either {\bf C/Z} or to
 a subgroup of {\bf R/Z}.}

\bigskip

\centerline{\bf References}

\bigskip

\noindent [Ba] I.N.BAKER, {\it Fractional iteration near a fixed
point of multiplier 1}, J. Austral. Math. Soc {\bf 4}, 1964, p.
143-148.

\medskip

\noindent [Ec] J.ECALLE, {\it Th\'eorie it\'erative: introduction
\`a la th\'eorie des invariantes holomorphes}, J. Math pures et
appliqu\'es {\bf 54}, 1975, p. 183-258.

\medskip

\noindent [Mo] J.MOSER, {\it On commuting circle mappings and
simultaneous Diophantine approximations}, Math. Zeitschrift {\bf
205}, 1990, p. 105-121.

\medskip

\noindent [PM1] R.PEREZ-MARCO, {\it Solution compl\`ete au
probl\`eme de Siegel de lin\'earisation d'une application
holomorphe au voisinage d'un point fixe (d'apr\`es J.C.Yoccoz)},
Seminar Bourbaki {\bf 753}, 1991.

\medskip

\noindent [PM2] R.PEREZ-MARCO, {\it Uncountable number of symmetries for non-linearisable holomorphic dynamics},
Inventiones Mathematicae, {\bf 119}, {\bf 1}, p.67-127,1995; (C.R. Acad. Sci. Paris, {\bf 313}, 1991, p.
461-464).

\medskip

\noindent [Vo] S.M.VORONIN, {\it Analytic classification of germs
of conformal mappings ({\bf C},0) $\to$ ({\bf C},0) with identity
linear part}, Funktsional Anal. i Prilozhen, {\bf 16}, 1982, p.
1-17.

\medskip

\noindent [Yo] J.C.YOCCOZ, {\it Petits diviseurs en dimension 1}, S.M.F., Asterisque {\bf 231} (1995)

\medskip
\medskip

\noindent Ramakrishna Mission Vivekananda University,
Belur Math, WB-711202, India.

\noindent E-mail: kingshook@rkmvu.ac.in

\end